\documentclass[12pt]{article}
\usepackage{amsmath}
\usepackage{graphicx}
\usepackage{amsfonts}
\usepackage{amssymb}
\usepackage{stmaryrd}
\usepackage{amscd}
\newcommand{\bc}{\begin{center}}
\newcommand{\ec}{\end{center}}
\newcommand{\be}{\begin{equation}}
\newcommand{\ee}{\end{equation}}
\newcommand{\bea}{\begin{eqnarray}}
\newcommand{\eea}{\end{eqnarray}}
\newcommand{\ba}{\begin{array}}
\newcommand{\ea}{\end{array}}
\newcommand{\edc}{\end{document}}

\def\l{\lambda}

\begin{document}
\thispagestyle{empty}
\begin{center}

{\textbf{ON DYNAMICS OF $\ell$- VOLTERRA QUADRATIC STOCHASTIC
OPERATORS}}\\[4mm]
{\bf U. A. Rozikov}$^{1,2}$ and {\bf A. Zada}$^2$\\[2mm]
$^1$Institute of Mathematics and Information Technologies,\\ Tashkent, Uzbekistan.\\
email: rozikovu@yandex.ru \\
$^2$Abdus Salam School of Mathematical Sciences, GCU,\\
35-C-2, Gulberg III, Lahore, Pakistan.\\
email: zadababo@yahoo.com\\
\vspace{0.5cm}
\end{center}
{\bf Abstract.} We introduce a notion of $\ell$-Volterra quadratic
stochastic operator defined on $(m-1)$-dimensional simplex, where
$\ell\in\{0,1,...,m\}$. The $\ell$-Volterra operator is a Volterra
operator iff $\ell=m$. We study structure of the set of all
$\ell$-Volterra operators and describe their several fixed and
periodic points. For $m=2$ and $3$ we describe behavior of
trajectories of $(m-1)$-Volterra operators. The paper also contains
many remarks with comparisons of  $\ell$-Volterra operators and
Volterra ones. {\vskip .5cm \noindent


\section{Introduction}
In biology a quadratic stochastic operator (QSO) has meaning of a
population evolution operator (see \cite{HS}, \cite{K}, \cite{L}),
which arises as follows. Consider a population consisting of $m$
species. Let $x^{0} = (x_{1}^{0},...,x_{m}^{0})$ be the probability
distribution of species in the initial generations, and $P_{ij,k}$
the probability that individuals in the $i$th and $j$th species
interbreed to produce an individual $k$. Then the probability
distribution $x'= (x_{1}',...,x_{m}')$ (the state) of the species in
the first generation can be found by the total probability i.e.
$$ x'_k=\mathop {\sum} \limits^{m}_{i,j=1}P_{ij,k}x^{0}_{i}x^{0}_{j} ,\,\,\,\,k=
1,...,m. \eqno(1)$$
 This means that the association $x^{0}\rightarrow x'$
defines a map $V$ called the evolution operator. The population
evolves by starting from an arbitrary state $x^{0}$, then passing to
the state $x'= V(x)$ (in the next "generation"), then to the state
$x''=V(V(x))$, and so on.
 Thus states of the population described by the following dynamical
 system
         $$ x^{0},\ \ x'= V(x), \ \ x''=V^{2}(x),\ \  x'''= V^{3}(x),...$$
Note that $V$ (defined by (1)) is a non linear (quadratic) operator,
and it is higher dimensional if $m\geq 3$. Higher dimensional
dynamical systems are important but there are relatively few
dynamical phenomena that are currently understood (\cite{De},
\cite{El}, \cite{R}).

In this paper we consider a class of nonlinear (quadratic) operators
which we call $\ell$-Volterra operators and study dynamical systems
generated by such operators.

The paper is organized as follows.

 In section 2 we give some preliminary definitions.
 Also we discuss the difference of
 quadratic operators introduced in this paper  from known quadratic operators.
 In section 3 we describe some invariant
 (in particular some fixed points) sets for $\ell$-Volterra operators.
 Also we give some family of $\ell$-Volterra operators each element
 of which has cyclic orbits generated by several
 vertices of the simplex. We also show that the set of all
 $\ell$-Volterra operators is convex,
 compact and describe its extremal points.
 Section 4 devoted to 1-Volterra operators and section 5 devoted to 2-Volterra
 operators defined on a two dimensional simplex.
In these sections we describe limit behavior of all trajectories
(orbits).

\section{Preliminaries}

The quadratic stochastic operator (QSO) is a mapping of the simplex.

$$ S^{m-1}=\left\{x=(x_1,...,x_m)\in {\bf R}^m: x_i\geq
0,\,\sum^m_{i=1}x_i=1 \right\} \eqno(2)$$ into itself , of the form
$$
    V: x'_k= \sum^m_{i,j=1}P_{ij,k}x_ix_j,\,\,\,k=
    1,...,m,\eqno(3)$$
where $P_{ij,k}$ are coefficients of heredity and
$$P_{ij,k}\geq 0,\ \ P_{ij,k}= P_{ji,k},\ \ \sum^m_{k=1}P_{ij,k}= 1, (i,j,k=
1,....,m),\eqno(4)$$

 Thus each quadratic stochastic operator $V$ can
be uniquely defined by a cubic matrix
$\textbf{P}=\left(P_{ij,k}\right)^n _{i,j,k=1}$ with conditions (4).

Note that each element $x\in S^{m-1}$ is a probability distribution
on $E= \{1,...,m\}$. The population evolves by starting from an
arbitrary state (probability distribution on $E$) $x\in S^{m-1}$
then passing to the state $V(x)$ (in the next "generation"), then to
the state $V(V(x))= V^2(x)$, and so on.

For a given $x^{(0)}\in S^{m-1}$ the trajectory (orbit)
$$ \{ x^{(n) }\}, \ \ n=0,1,2,...\ \ \mbox{of}\ \ x^{(0)}$$
under the action of QSO (3) is defined by
$$x^{(n+1)}= V(x^{(n)}),\ \ \mbox{where} \ \ n= 0,1,2,...$$

One of the main problem in mathematical biology consists in the
study of
 the asymptotical behavior of the trajectories. The difficulty of the problem
 depends on given matrix \textbf{P}. Now we shall briefly describe the history of (particularly)
 studied QSOs, which allows the reader to easily see the place of the operators introduced in this paper.

{\it The Volterra operators.}(see \cite{Gan1},
\cite{Gan3},\cite{GE}) A Volterra QSO is defined by (3), (4)
 and the additional assumption
$$ P_{ij,k}=0, \ \ \mbox{if}\ \ k\not\in \{i,j\},\,\,
\forall i,j,k\in E. \eqno(5)$$

The biological treatment of condition (5) is clear: The offspring
repeats the genotype of one of its parents.

In paper \cite{Gan1} the general form of Volterra QSO
$$ V: x= (x_{1},...,x_{m}) \in S^{m-1}\,
\rightarrow\,V(x)= x'= (x'_1,...,x'_m)\,\in S^{m-1}
$$ is given
$$x'_k=x_k\left(1+\sum^m_{i=1}a_{ki}x_i\right),\eqno(6)$$
where
$$a_{ki}=2P_{ik,k}-1 \ \ \mbox{for}\, i\neq k\,\,\mbox{and}\
\,a_{ii}=0, i=1,...,m.$$  Moreover
$$a_{ki}=-a_{ik}\ \ \mbox{and} \ \ |a_{ki}| \leq 1.$$

  In \cite{Gan1}, \cite{Gan3} the theory of QSO (6) was developed
by using theory of the Lyapunov function and tournaments. But
non-Volterra QSOs (i.e. which do not satisfy the condition (5)) were
not in completely studied. Because there is no any general theory
which can be applied for investigation of non-Volterra operators.

 To the best of our knowledge, there are few papers devoted to such
 operators. Now we shall describe non-Volterra  operators:

{\it A permutated Volterra operator.} Papers \cite{GD}, \cite{GA}
are devoted to study of non-Voltera operators which are generated
from Volterra operators (6) by a cyclic  permutation of coordinates
i.e $$ V_{\pi} : x'_{\pi(i)}=x_{i}\left(1+\sum^m_{i=1}a_{ki}
x_i\right), \ \ i= 1,...,m,$$ where  $\pi$ is a cyclic permutation
on the set of indices $E$.

{\it Kvazi-Volterra QSO.} In \cite{GM} a class of Kvazi-Volterra
operators is introduced. For such operators the condition (5) is not
satisfied only for very few values of $i,j,k$.

{\it Non-Volterra QSO as combination of a Volterra and a
non-Volterra operators.} In paper \cite{Gan} it was considered the
following family of QSOs
$V_{\lambda}:S^{2}\,\,\rightarrow\,\,S^{2}:$
$$V_{\lambda}=(1-\lambda)V_0+\lambda V_1, \,\,0\leq\lambda\leq 1,$$
 where
$$V_0(x)=\left(x^2_1+2x_1x_2,\,x^2_2+2x_2x_3,\,x^2_3+2x_1x_3\right),$$
is Volterra operator and
$$ V_1(x)=\left(x^{2}_{1} + 2x_{2}x_{3},\,x^{2}_{2} + 2x_{1}x_{3},\,x^{2}_{3} +
2x_{1}x_{2}\right),$$ is non-Volterra QSO.

Note that behavior of the trajectories of $V_{0}$ is very irregular
(see \cite{L}, \cite{SU}, \cite{Z}).

{\it Non-Volterra QSO generated by a product measure.} In \cite{Ga},
\cite{GR} a constructive description of the matrix \textbf{P} is
given. This construction depends on a probability measure $\mu$
which is given on a fixed graph $G$ and cardinality of a set of
cells (configurations).

 In \cite{Ga} it was proven that the QSO constructed by
the construction is Volterra iff $G$ is a connected graph.

In \cite{RS} using the construction of QSO for the general finite
graph and probability measure $\mu$ (here $\mu $ is product of
measures defined on maximal subgraphs of the graph $G$) a class of
non-Volterra QSOs is described.

 It was shown that if $\mu$ is given by the
product of the probability measures then corresponding non-Volterra
operators can be studied by $N$ number (where $N$ is the number of
maximal connected subgraphs) of Volterra operators defined on the
maximal connected subgraphs.

{\it F-QSO.} Consider $E_0=E\bigcup \{0\}=\{0,1,...,m\}.$ Fix a set
$F\subset E$ and call this set the set of "females" and the set
$M=E\setminus F$ is called the set of "males".  The element $0$ will
play the role of "empty-body". Coefficients $P_{ij,k}$ of the matrix
${\mathbf P}$ we define as follows
$$P_{ij,k}=\left\{\begin{array}{lll}
1, \ \ {\rm if} \ \ k=0, i,j\in F\cup \{0\} \ \ {\rm or} \ \ i,j\in
M\cup \{0\};\\
0, \ \ {\rm if} \ \ k\ne 0, i,j\in F\cup \{0\} \ \ {\rm or} \ \
i,j\in
M\cup \{0\};\\
\geq 0, \ \ {\rm if } \ \ i\in F, j\in M, \forall k.
\end{array}\right.\eqno (7)$$

Biological treatment of the coefficients (7) is very clear: a \ \
"child" $k$ can be generated if its parents are taken from different
classes $F$ and $M$.

For a given $F\subset E$ a QSO with (3),(4) and (7) is called a
$F$-QSO. Note that $F$-QSO is non-Volterra for any $ F \subset E$.

 In \cite{RJ} the $F$-QSOs are studied for any $F\subset E$.
 It was proven that such operators has unique fixed
point $(1,0,...,0) \in S^{m}$ and all trajectories converge to this
fixed point faster than any geometric progression.

{\it Strictly non-Volterra QSO.} Recently in \cite{RJ1} a new class
of non-Volterra operators is introduced. These operators satisfy
$$ P_{ij,k}=0,\ \ \mbox{if} \ \ k\in \{i,j\}, \ \ \forall i,j,k\in
E. \eqno(8)$$ Such an operator is called strictly non-Volterra QSO.

For arbitrary strictly non-Volterra QSO defined on $S^{2}$ in
\cite{RJ1} it was proved that every such an operator has a unique
fixed point. Also it was proven that, such operators have a cyclic
trajectory. This is quite different behavior from the behavior of
Volterra operators, since the Volterra operators have no any cyclic
trajectory.

Now we shall give a new class of non-Volterra operators.

{\it $\ell$-Volterra QSO.} Fix $\ell \in \{1,...,m\}$ and assume
that elements $P_{ij,k}$ of the matrix $\textbf{P}$ satisfy
$$ P_{ij,k}= 0 \ \ \mbox{if} \ \ k \not\in
\{i,j\}\ \ \ \mbox{for}\ \ \mbox{any}\ \ k\in \{1,...,\ell\},\ \ i,j
\in E;\eqno(9)$$
$$ P_{ij,k}> 0 \ \ \mbox{for at least one pair} \ \ (i,j),\ \
i\neq k,\ \ j\neq k \ \ \mbox{if} \ \  k \in
\{\ell+1,...,m\}.\eqno(9a)$$

\textbf{Definition 1.} {\it For any fixed $\ell \in \{1,...,m\}$,
the QSO defined by (3), (4), (9) and (9a) is called $\ell$-Volterra
QSO.}

Denote by $\mathcal{V}_{\ell}$ the set of all $\ell$-Volterra QSOs.

\textbf{Remarks.} 1. The condition (9a) guarantees that $\mathcal{V}
_{\ell_{1}} \bigcap \mathcal{V}_{\ell_{2}} = \emptyset $ for any
$\ell_{1}\neq \ell _{2}$.

2. Note that $\ell$-Volterra QSO is Volterra if and only if $\ell=
m$.

3. Kvazi-Volterra operators (introduce above) are particular case of
$\ell$-Volterra operators.

4. The class of $\ell$-Volterra QSO for a given $\ell$ does not
coincide with a class of non-Volterra QSOs mentioned above.

We shall use the following notions.
 \vskip 0.3 truecm

{\bf Definition 2.}(\cite{De}, p. 215 ){\it A fixed point $P$ for $
F: {\mathbf R}^m \rightarrow {\mathbf R}^m$ is called hyperbolic if
the Jacobian matrix $\textbf{J}$ of the map $F$ at the point $P$ has
no eigenvalues on the unit circle.}

There are three types of hyperbolic fixed points :

1. $P$  is an {\it attracting} fixed point if all of the eigenvalues
of $\textbf{J}(P)$ are less than one in absolute value.

2. $P$ is a {\it repelling} fixed point if all of the eigenvalues of
$\textbf{J}(P)$ are greater than one in absolute value.

3. $P$ is a {\it saddle} point otherwise.

The following theorem is also very useful.
 \vskip 0.3 truecm

{\bf Theorem 1.} (\cite{De}, p.217) {\it Suppose $F$ has a saddle
fixed point $P$. There  exist $\varepsilon > 0$ and a smooth  curve
$\gamma :(-\varepsilon ,\varepsilon) \rightarrow {\mathbf R}^2$ such
that $\gamma(0) = P$; $\gamma'(t)\neq 0$; $\gamma'(0)$ is an
unstable eigenvector for $\textbf{J}(P)$; $\gamma $ is $F^{-1}-$
invariant; $F^{-n}(\gamma(t)) \rightarrow P$ as $n \rightarrow
\infty $; if $|F^{-n}(q)- P|< \varepsilon$ for all $n \geq 0$ then
$q=\gamma(t)$ for some $t$.}

 The curve $\gamma$ is called the (local) unstable manifold at $P$.
 The theorem is true for stable sets as well as with the obvious
 modification. On the local manifold, all points tend to the fixed
 point under iteration of $F$.

\section{Canonical form of $\ell$-Volterra QSO.}

By definition  for  $k=1,...,\ell$ we have
 $$x'_k=\sum^m_{i,j=1}P_{ij,k}x_ix_j=P_{kk,k}x^2_k+2\sum^m_{i=1\atop i\neq
 k}P_{ik,k}x_ix_k=$$
$$x_k\left(P_{kk,k}x_k + 2\sum^m_{i=1\atop i\neq
k}P_{ik,k}x_i\right).$$

 Using  $x_k=1-\sum^m_{i=1\atop i\neq k}x_{i}$ we get
 $$x'_k =x_k\left(P_{kk,k}+\sum^m_{i=1\atop i\neq k}\left(2P_{ik,k}-P_{kk,k}\right)x_{i}\right),\ \
k=1,...,\ell.$$ For $k= \ell+1,...,m$ we have
$$ x'_k =x_k\left(P_{kk,k}+\sum^m_{i= 1\atop i\neq k}\left(2P_{ik,k}-P_{kk,k}\right)x_i\right)
+\sum^m_{{i,j= 1\atop i\neq k}\atop j\neq k}P_{ij,k}x_ix_j.$$ Denote
$a_{ki}=2P_{ik,k}-P_{kk,k}$ then

$$\left\{\begin{array}{llllll}
x_k' =x_k\left(a_{kk}+\sum^m_{i=1\atop i\neq k}a_{ki}x_i\right),\ \
k=
1,...,\ell\\[3mm]
x_k' =x_k\left(a_{kk}+\sum^m_{i=1\atop i\neq
k}a_{ki}x_i\right)+\sum^m_{{i,j= 1\atop i\neq k}\atop j\neq
k}P_{ij,k}x_ix_j, \ \ k= \ell+1,...,m.
\end{array}\right.\eqno (10)$$

Note that $0\leq a_{kk} \leq 1$ and $-a_{kk}\leq a_{ki} \leq
2-a_{kk},$ $i\neq k$, $0\leq P_{ij,k} \leq 1.$

For any $I \subset E= \{1,2,...,m\}$ we define the face of the
simplex $S^{m-1}$:
$$\Gamma_I= \left\{x \in S^{m-1}: x_i= 0 \ \
\mbox{for any} \ \ i\in I \right\}.$$ \vskip 0.3 truecm

{\bf Proposition 1.} {\it Let $V$ be a $\ell$-Volterra QSO. Then the
following are true

(i) Any face $\Gamma_{I}$ with $I\subseteq \{1,...,\ell \}$ is
invariant set with respect to $V.$

(ii) Let $A_{\ell} = \{i \in \{1,..., \ell \}: a_{ii}> 0\}$. For any
$I\subset A_{\ell} \cup \,\{ \ell + 1 ,...,m\}$ the set $T_{I}= \{ x
\in S^{m - 1 }: x_{i} > 0, \forall i \in I\}$ is invariant with
respect to $V$.}
 \vskip 0.3truecm

 {\bf Proof.} (i) From (10) it follows that if $x_i= 0$ then $x'_i = 0$ for
any $i \in \{1,...,\ell\}$. Hence $V(\Gamma_{I}) \subset \Gamma_{I}$
if  $I \subset \{1,...,\ell \}$.

(ii) Take $I \subset A_{\ell} \cup \{\ell+1,...,m\}$. For $k \in I
\cap A_{\ell}$ by (10) and inequality $-a_{kk} \leq a_{kj}, j= 1
,...,m$ we get
$$ x'_k = x_k\left(a_{kk}+\sum^m_{j=1\atop j\ne
k}a_{kj}x_j\right)\geq$$ $$ x_k\left(a_{kk}-a_{kk}\sum^m_{j=1\atop
j\ne k}x_j\right)=a_{kk}x^2_k>0,\ \ \mbox{since} \ \ x_k>0\ \
\mbox{for}\ \ k \in I \cap A_{\ell}.$$

For $k\in I \cap\{\ell+1,...,m\}$ by (10) and  condition (9a) we
have
$$ x'_k = x_k\left(a_{kk}+\sum^m_{j=1\atop j\ne k}a_{kj}x_j\right)+
\sum^m_{{i,j=1\atop i\neq k}\atop j\neq k}P_{ij,k}x_ix_j\geq$$
$$
x_k\left(a_{kk}-a_{kk}\sum^m_{j= 1\atop j\ne k}x_j\right)
+\sum_{{i,j\in I\atop i\neq k}\atop j\neq k}P_{ij,k}x_ix_j
> a_{kk}x^2_k \geq 0,$$
here we used $\sum_{{i,j\in I\atop i\neq k}\atop j\neq
k}P_{ij,k}x_ix_j> 0$ which follows from condition (9a) and $x_i> 0,
x_j>0, \ \ \forall i,j\in I$. Thus $V(T_{I})\subset T_{I}$ if $I
\subset A_{\ell}\cup \{\ell+1,...,m\}$. The proposition is proved.
\vskip 0.3 truecm

 Denote $e_i=\left(\delta_{1i},...,\delta_{mi}\right)\in S^{m-1}$,
 $i=1,...,m$ the vertices of the simplex $S^{m-1}$,
 where $\delta_{ij}$ is the Kronecker's symbol.
 \vskip 0.3 truecm

{\bf Proposition 2.} 1) {\it The vertex $e_{i}$ is a fixed point for
a $\ell$-Volterra QSO iff} $P_{ii,i}=1$, $(i= 1,...,m)$.

2) {\it For any collection $I_s=\{e_{i_{1}},...,{e_{i_{s}}}\}
\subset \{e_{\ell + 1},...,{e_{m}}\},\ \ (s\leq m - \ell )$ there
exists a family $\mathcal{V}_{\ell}(I_s) \subset \mathcal{V}_{\ell}$
such that
 $\{e_{i_{1}},...,{e_{i_{s}}}\}$ is a $s$-cycle for each} $V\in \mathcal{V}_{\ell}(I_s)$.
\vskip 0.3 truecm

{\bf Proof.} 1) It is easy to see that if $i\in \{1,...,\ell\}$ then
$$ V(e_i)=\left(0,...,0,P_{ii,i},0,...,0,P_{ii,\ell+1},...,P_{ii,m}\right)\ \
\hbox{with} \ \ P_{ii,i}+\sum^m_{j=\ell+1}P_{ii,j}=1$$ and if $i\in
\{\ell+1,...,m\}$ then
$$V(e_i)=\left(0,...,0,P_{ii,\ell+1},...,P_{ii,m}\right) \ \
\mbox{with} \ \ \sum^m_{j=\ell+1}P_{ii,j} = 1.\eqno (11)$$ Thus
$V(e_i)=e_i$ iff $P_{ii,i}=1.$

2) By (11) we have
$$V(e_{i_{j}})=\left(0,...,0,P_{i_{j}i_{j},\ell+1},...,P_{i_{j}i_{j},m}\right)$$ for any $j=
 1,...,s$. In order to get $V(e_{i_{1}})=e_{i_{2}}$ we assume
 $$P_{i_{1}i_{1},i_{2}}=1, \ \ P_{i_{1}i_{1},j}= 0, \ \ j\neq i_{2}. \eqno(11_1)$$
 Then to get $V(e_{i_{2}})=e_{i_{3}}$ we assume
 $$P_{i_{2}i_{2},i_{3}}=1, \ \ P_{i_{2}i_{2},j}=0, \ \ j\neq i_{3}. \eqno(11_{2})$$

 Similarly to get $V(e_{i_{s-1}})= e_{i_{s}}$ we assume
 $$P_{i_{s-1}i_{s-1},i_{s}}=1, \ \ P_{i_{s-1}i_{s-1},j}=0, \ \ j\neq i_{s}. \eqno(11_{s-1})$$
 The last assumption follows from $V(e_{i_{s}})=e_{i_{1}}$ i.e
 $$ P_{i_{s}i_{s},i_{1}}=1, \ \ P_{i_{s}i_{s},j}= 0,\ \ j\neq i_{1}. \eqno(11_s)$$

 Hence $\mathcal{V}_{\ell}(I_s)= \{V \in  \mathcal{V}_{\ell}: \hbox{ the coefficients of $V$
 satisfy} \ \ (11_1)-(11_s)\}$.
 The proposition is proved.

For any set $A$ denote by $|A|$ its cardinality.

 The next proposition gives
 a set of periodic orbits of $\ell$-Volterra QSOs.
\vskip 0.3 truecm

\textbf{Proposition 3.} {\it For any $I_1,...,I_q\subset \{\ell+1
 ,...,m\}$ such that $I_{i} \bigcap I_{j}= \emptyset \ \ (i\neq j, i,j=
 1,...,q)$. There exists a family $\mathcal{V}_{\ell}(I_1,...,I_q)\subset \mathcal{V}_{\ell}$
such that each collection $\{e_i, i\in I_j\}, \  \ j= 1,...,q$ is a
$|I_j|-$ cycle for every $V \in \mathcal{V}_{\ell}(I_1,...,I_q)$.}
\vskip 0.3 truecm

{\bf Proof.} Since $I_i \bigcap I_j= \emptyset$, $i\neq j$ the
 family can be constructed using Proposition 2 i.e. $\mathcal{V}_{\ell}(I_1,...,I_q)=
 \bigcap^q_{i=1}\mathcal{V}_{\ell}(I_i)$.
\vskip 0.3 truecm

\textbf{Remarks.} 1) There is not any $\ell $-Volterra operator with
a periodic orbit  $\{e_{i_{1}},...,{e_{i_{s}}}\} \subset
\{e_{1},...,{e_{\ell }}\},\,1< s\leq \ell
 $.

 2) Propositions 2 and 3 show that $\ell$-Volterra operators
 have quite different behavior from the behavior of Volterra operators,
 since Volterra operators have no cyclic trajectories.

Recall that $\mathcal{V}_{\ell}$ is the set of all $\ell$-Volterra
operators defined on $S^{m - 1}$. \vskip 0.3 truecm

\textbf{Proposition 4.} (i) {\it The set  $\mathcal{V}_{\ell}$ is a
convex, compact subset of} $\textbf{R}^\frac{m(m-1)(m-\ell +
1)}{2}.$

(ii) {\it The extremal points of $\mathcal{V}_{\ell}$ are
$\ell$-Volterra operators with $P_{ij,k}= 0$ or $1$ for any $i,j,k$
i.e.}
$${\rm Extr}(\mathcal{V}_{\ell})=\left\{V\in \mathcal{V}_{\ell}:
\mbox{the matrix} \ \ \mathbf{P} \ \ \mbox{of} \ \ V \ \
\mbox{contains only} \ \ 0 \ \ \mbox{and} \ \ 1\right\}.$$

(iii) {\it If $\ell=m$ then $\left|{\rm
Extr}(\mathcal{V}_{\ell})\right|=2^{\frac{1}{2}m(m-1)}$; if
$\ell\leq m-1$ then}
$$\left|{\rm Extr}(\mathcal{V}_{\ell})\right|=\big(m-\ell\big)^{\frac{1}{2}(m-\ell)(m-\ell+1)}
\big(m-\ell+1\big)^{(m-\ell+1)\ell}\big(m-\ell+2\big)^{\frac{1}{2}\ell(\ell-1)}.$$
 \vskip 0.3 truecm

 {\bf Proof.} (i) Since we have one-to-one correspondence between the set of all
 QSOs and the set of all cubic matrices $\mathbf{P}$,
 we can consider a QSO $V$ as a point of $\mathbf{R}^{m(m^2-1)}.$
 The number $\frac{m(m-1)(m-\ell + 1)}{2}$ is the number of
 independent elements of the matrix $\textbf{P}$ with the condition
 (9). Let $V_1, V_2$ be two $\ell$-Volterra QSO i.e $V_{1},V_{2}\in \mathcal{V}_{\ell}$.
 We shall prove that $V=\lambda V_{1} + (1-\lambda) V_{2}\in \mathcal{V}_{\ell}$ for
 any $\lambda \in [0,1]$.

 Let $P^{(1)}_{ij,k}$ (resp. $P^{(2)}_{ij,k})$ be coefficients of $V_{1}$ (resp. $V_{2})$.
 Then coefficients of $V$ has the form
$$P_{ij,k}= \lambda P^{(1)}_{ij,k} + ( 1 - \lambda )P^{(2)}_{ij,k}. \eqno(12)$$
By definition coefficients $ P^{(1)}_{ij,k} $ and $P^{(2)}_{ij,k}$
satisfy conditions (9) and (9a). Using $(12)$ it is easy to check
that $P_{ij,k}$ also satisfy the condition (9) and (9a).

(ii) Assume $V \in \mathcal{V}_{\ell}$ with $ P_{i_{0}j_{0},k_{0}}=
\alpha \neq 0$ and $1$ for some $i_{0},j_{0},k_{0}$. Construct two
operators $V_{q}$ with coefficients $ P^{(q)}_{ij,k}$, $q= 1,2$ as
following
$$ P^{(1)}_{ij,k}=  \left\{\begin{array}{ll}
P_{ij,k} \ \  \mbox{if} \ \ {(i,j)}\neq {(i_{0},j_{0})}, \\
1 \ \ \ \ \mbox{if} \ \ {(i,j,k)}= {(i_{0},j_{0},k_{0})},\\
0 \ \ \ \ \mbox{if} \ \ {(i,j,k)}= {(i_{0},j_{0},k)}, \ \ k\neq
k_{0}
\end{array}\right. $$

$$ P^{(2)}_{ij,k}=  \left\{\begin{array}{ll}
P_{ij,k} \ \ \mbox{if} \ \ {(i,j)}\neq {(i_{0},j_{0})}, \\
0 \ \ \ \ \mbox{if} \ \ {(i,j,k)}= {(i_{0},j_{0},k_{0})},\\
\frac{P_{ij,k}}{1 - \alpha} \ \ \mbox{if} \ \ {(i,j,k)}=
{(i_{0},j_{0},k)}, \ \ k\neq k_{0}.
\end{array}\right.$$
Then
$$ \alpha P^{(1)}_{ij,k} + ( 1 - \alpha ) P^{(2)}_{ij,k}=  \left\{\begin{array}{ll}
P_{ij,k} \ \ \mbox{if} \ \ {(i,j)}\neq {(i_{0},j_{0})}, \\
\alpha =P_{i_{0}j_{0},k_{0}} \ \ \mbox{if} \ \ {(i,j,k)}=
{(i_{0},j_{0},k_{0})}\hspace {0.6cm} =\ \ P_{ij,k}.\\
P_{ij,k} \ \ \ \ \mbox{if} \ \ {(i,j,k)}= {(i_{0},j_{0},k)}, \ \
k\neq k_{0}
\end{array}\right.\eqno(12')$$
Since $\alpha>0$, from (12') we get $P_{ij,k}=0$ if and only if
$P^{(1)}_{ij,k}=0$ and $P^{(2)}_{ij,k}=0.$ This means that $V_1$ and
$V_2$ are $\ell$-Volterra operators. Hence $V= \alpha V_{1} + ( 1 -
\alpha) V_{2}$. Thus if $P_{ij,k}\in (0,1)$ for some $(i,j,k)$ then
$V$ is not an extremal point. Finally, if $P_{ij,k}= 0$ or 1 for any
$(i,j,k)$ then the representation $V=\lambda V_1+(1-\lambda)V_2$,
$0<\lambda<1$ is possible only if $V_1= V_2= V$.

(iii) In order to compute cardinality of ${\rm
Extr}(\mathcal{V}_\ell)$ we have to know which elements of the
matrix $\mathbf{P}$ can be 1.

Denote $\mathcal{P}_{ij}=\left(P_{ij,1},...,P_{ij,m}\right)^t$ the
$(i,j)$th column of $\mathbf{P}$, where $(i,j)\in
\mathcal{K}=\{(i,j): 1\leq i\leq j\leq m\}.$

Let $n_0(\mathcal{P}_{ij})$ be the number of elements of
$\mathcal{P}_{ij}$ which must be zero by conditions (4), (9), (9a).

Put for $\ell\in \{1,...,m\}:$

$$\mathcal{A}\equiv \mathcal{A}_{em}=\left\{(i,j)\in \mathcal{K}: i\leq
\ell, j\in \{i\}\cup \{\ell+1,...,m\}\right\},$$
$$\mathcal{B}\equiv \mathcal{B}_{em}=\left\{(i,j)\in \mathcal{K}: i\leq
\ell, j\leq \ell, i<j\right\},$$
$$\mathcal{C}\equiv \mathcal{C}_{em}=\left\{(i,j)\in \mathcal{K}: \ell<i\leq j\right\}.$$
Note that  $\mathcal{K}=\mathcal{A}\cup\mathcal{B}\cup\mathcal{C}.$
If $\ell=m$ then $\mathcal{C}=\emptyset.$

It is easy to see that
$$n_0(\mathcal{P}_{ij})=\left\{\begin{array}{lll}
\ell-1 \ \ \mbox{if} \ \ (i,j)\in \mathcal{A}\\
\ell-2 \ \ \mbox{if} \ \ (i,j)\in \mathcal{B}\\
\ell \ \ \ \ \mbox{if} \ \ (i,j)\in \mathcal{C}\\
\end{array}\right.$$

By condition (4) each column contains unique "1". We have
$m-n_0(\mathcal{P}_{ij})$ possibilities to write 1 in the column
$(i,j)\in \mathcal{K}.$ Thus
$$\left|{\rm Extr}(\mathcal{V}_\ell)\right|=\prod_{(i,j)\in
\mathcal{K}}\left(m-n_0(\mathcal{P}_{ij})\right)=(m-\ell+1)^{|\mathcal{A}|}(m-\ell+2)^{|\mathcal{B}|}
(m-\ell)^{|\mathcal{C}|}.$$ This with
$$|\mathcal{A}|=(m-\ell+1)\ell,\ \ |\mathcal{B}|=\frac{1}{2}(\ell-1)\ell,\ \
|\mathcal{C}|=\frac{1}{2}(m-\ell+1)(m-\ell),$$ would yield the
formula. The proposition is proved.

For the set $\mathcal{V}$ of all QSOs we have  $\mathcal{V} \subset
\textbf{R}^{\frac{m{(m^{2} - 1)}}{2}}$. Note that $\mathcal{V}$ also
is a convex, compact set. Its extremal points also are operators
with $P_{ij,k}=0$ or 1 only. It is easy to see that

$$\left|{\rm Extr}(\mathcal{V}_m)\right|<\left|{\rm
Extr}(\mathcal{V}_{m-1})\right|<...<\left|{\rm
Extr}(\mathcal{V}_1)\right|<\left|{\rm
Extr}(\mathcal{V})\right|=m^{\frac{1}{2}m(m+1)}.$$
 For example, if $m=3$ then
 $$\left|{\rm Extr}(\mathcal{V}_3)\right|=8, \ \ \left|{\rm
 Extr}(\mathcal{V}_2)\right|=48, \ \ \left|{\rm
 Extr}(\mathcal{V}_1)\right|=216, \ \ \left|{\rm
 Extr}(\mathcal{V})\right|=729.$$

The set $\mathcal{V}$ can be written as $\mathcal{V}= \bigcup
^m_{\ell =0} \mathcal{V}_{\ell}$. Here $\mathcal{V}_{0}$ is the set
of " $0$-Volterra QSO"s i.e for any $k\in \{1,...,m\}$ there is at
least one pair $(i,j)$ with $i\neq k$ and $j\neq k$ such that
$P_{ij,k}>0 $.

As it was mentioned above:  $\mathcal{V}_m$ is the set of all
Volterra operators and $\mathcal{V}_{\ell_1} \bigcap
\mathcal{V}_{\ell_2} = \emptyset$ for any $\ell_1 \neq \ell_2 \in
\{0,...,m\}.$

Thus to study dynamics of QSOs from $\mathcal{V}$ it is enough to
study the problem for each $\mathcal{V}_{\ell}$ , $\ell = 0,...,m. $

In general, the problem of study the behavior of $V \subset
\mathcal{V}_{\ell}$ (for fixed ${\ell}$) is also a difficult
problem. So in the next sections we consider the problem for small
dimensions (i.e $m= 2,3$) and $\ell =1,2$. \vskip 0.3 truecm

\section{Case $m= 2$} In the case $m=2$ we have only
$1$-Volterra operator $V : S^1 \rightarrow S^1$ such that
$$ \left\{\begin{array}{ll}
x'= ax^{2} + 2cxy\\
y'= bx^{2} + 2dxy + y^{2},\\
\end{array}\right.\eqno(13)$$
where $a,b,c,d \in[0,1)$ (the case $a=1$ corresponds to Volterra
operator), $ a + b = c + d = 1$. Using $x + y = 1$ from (13) we get
a dynamical system generated by function $f(x) = (a - 2c)x^{2} +
2cx, x\in [0,1], a \in [0,1), c \in [0,1]$. By properties of $f(x)$
one can prove the following
 \vskip 0.3 truecm

 {\bf Proposition 5.} 1) {\it If $c\leq \frac{1}{2}$ , $\forall a \in
 [0,1)$ the operator (13) has unique fixed point $\lambda_{0}
= (0,1)$ and for any initial point $\lambda^{0} = (x^{0},y^{0}) \in
S^1$ the trajectory $\lambda^{(n)}$ goes to $\lambda_{0}$ as
$n\rightarrow \infty $.}

2) {\it If $c>\frac{1}{2}, \ \ \forall a \in
 [0,1)$ then (13) has two fixed points $\lambda_{0}
= (0,1)$ and $\lambda^{*} = (\frac{2c-1}{2c-a},\frac{1-a}{2c-a})$
the point $\lambda_{0}$ is repeller. For any initial point
$\lambda^0 \in S^1\setminus \{\lambda_{0}\}$ the trajectory
$\lambda^{(n)}$ tends to $\lambda^{*}$ as} $n\rightarrow\infty$.
\vskip 0.3 truecm

\section{Case $m=3$} In case $m= 3$ one has two
$\ell$-Volterra operators (for $\ell = 2$ and $1$). Here we shall
study the $2$-Volterra operators.

Arbitrary $2$-Volterra operator (for $m=3$) has the form :
$$ \left\{\begin{array}{ll}x'= x(a_{1}x + 2b _{1}y +
2c_{1}z)\\y'= y(2b_{2}x + d_{1}y + 2e_{1}z) \\z'= z(2c_{2}x +
2e_{2}y + z) + a_{2}x^{2} + 2b_{3}xy + d_{2}y_{2},
\end{array}\right.\eqno{(14)} $$
where $$a_{1}= P_{11,1}, \ \ a_{2}=P_{11,3};\ \ b_{i}=P_{12,i}, i=
1,2,3; \ \ c_{1}=P_{13,1},$$ $$ c_{2}=P_{13,3};\ \ d_{i}=P_{22,i},
i= 2,3; \ \ e_{i}= P_{23,i}, i= 2,3. \eqno(15)$$

To avoid many special cases and complicated formulas we consider the
case
$$
P_{11,1}=P_{22,2}, \ \ P_{13,1}=P_{23,2}, \ \
P_{12,1}=P_{12,2}.\eqno(16)$$ This corresponds to a symmetric (with
respect to permutations of $1$ and $2$) model.

Using $x + y + z= 1$ and condition (16) the operator (14) can be
written as
$$ \left\{\begin{array}{ll} x'= x (2c + (a-2c)x + 2(b-c)y)\\
y'= y(2c + 2(b-c)x + (a - 2c)y),
\end{array}\right.\eqno{(17)} $$
where $a= P_{11,1} \in [0,1),$ $b= P_{12,1} \in [0,\frac{1}{2}],$
$c= P_{13,1} \in [0,1],$ and $x,y \in [0,1]$ such that $x + y \leq
1$. \vskip 0.3 truecm

 {\bf Remark.} The case $a= P_{11,1}= P_{22,2}= 1$ corresponds to the
 Volterra case, so we consider only $a\neq 1$.
\vskip 0.3 truecm

 {\bf Theorem 2.} (i) {\it For $c\leq \frac{1}{2}$ the operator
 (17) has unique fixed point $\lambda _{0}= (0,0)$ which is
  global attractive point.}

 (ii) {\it Sets $M_{0} = \{\lambda =(x,y) : x= 0\}$, $M_{1} = \{\lambda=(x,y): y=
 0\}$, $M_{=} = \{\lambda =(x,y): x= y\}$, $M_{>} = \{\lambda =(x,y) : x> y\}$,
 $M_{<} =\{\lambda =(x,y): x< y\}$ are invariant with
 respect to the operator} (17).

(iii) {\it For $c>\frac{1}{2}$, $a\neq 2b$ the operator (17) has
four fixed points $\lambda _{0}=(0,0)$, $\lambda
_{1}=\left(0,\frac{2c-1}{2c-a}\right)$, $\lambda
_{2}=\left(\frac{2c-1}{2c-a},0\right)$,
$\lambda_3=\left(\frac{1-2c}{a + b - 4c},\frac{1 - 2c }{a + 2b -
4c}\right)$. Moreover $\lambda_0$ is repeller and
$$ \lambda_1 \ \ \hbox{and}\ \ \lambda_2 \ \ \mbox{are}\ \
\left\{\begin{array}{ll} \mbox{attractive, if} \ \  a>2b\\
\mbox{saddle, if} \ \ a<2b\\
\end{array}\right.$$

$$\lambda_3\ \ \hbox{is}\ \ \left\{\begin{array}{ll}
\mbox{attractive, if}\ \ a < 2b\\
\mbox{saddle, if} \ \ a > 2b.\\
\end{array}\right.$$}

 (iv) {\it For $c > \frac{1}{2}$ , $a = 2b$ the operator (17) has a
repeller fixed point $\lambda _{0} = (0,0)$ and continuum set of
fixed points $ F= \{ \lambda = (x,y) : x + y = \frac{2c - 1}{2(c-b)
}\}.$ The following line
$$I_\nu = \left\{\lambda=(x,y): y = \nu x, x\in [0,1] \right\}$$
is an invariant set for any $\nu \in [0,\infty)$. If $\lambda
^{0}=(x^0,y^0)$ is an initial point with $\frac{y^{0}}{x^{0}}= \nu$,
$(x^{0} \neq 0)$ then its trajectory $\lambda ^{(n)}$ goes to
$\overline{\lambda}
_{\nu}=\left(\frac{2c-1}{2(c-b)(1+\nu)},\frac{(2c-1)\nu}{2(c-b)(1
+\nu)}\right) \in I_{\nu} \bigcap F$ as $n\rightarrow \infty$, $\nu
\in [0,\infty)$, (if $x^{0}= 0$ then on invariant set $M_{0}$ we
have $\lambda ^{(n)}\rightarrow \lambda_1$).}

(v) {\it If $a<2b$ then $M_{0}$ (resp. $M_{1})$ is the stable
manifold of the saddle point $\lambda_{1}$ (resp. $\lambda_2)$. If
$a> 2b$ then $M_{=}$ is the stable manifold of saddle point
$\lambda_3$. There is an invariant curve $\gamma$ passing through
$\lambda_1,\lambda_2, \lambda_3$ which is unstable manifold for the
saddle points.}
 \vskip 0.3 truecm

 {\bf Proof.} (i) Clearly $\lambda _{0} = (0,0)$ is a fixed point
 for (17). Note that the Jacobian of (17) at (0,0) has the form
$$  \textbf{J}=\left(
\begin{array}{cc}
2c & 0 \\[2mm]
0 & 2c
\end{array}
\right), $$ so $\lambda_0$ is an attractive if $c< \frac{1}{2}$ and
non-hyperbolic if $c = \frac{1}{2}$.

Now we shall prove (for $c\leq \frac{1}{2}$) its global
attractiveness. From the first equation of (17) we have
$$x' = x( ax + 2by + 2cz ) \leq qx, \eqno(18)$$
where $q = \max\{a,2b,2c \}$. By definition of the operator (17) and
condition $c \leq \frac{1}{2}$ we have $q \leq 1$. Consider two
cases: \vskip 0.3 truecm

{\bf Case $q<1$.} In this case from (18) we get $x_{n + 1 }\leq
qx_{n} \leq q^{n}x^0$, where $x_{n}$ is the first coordinate of the
trajectory $\lambda^{(n)}=V^n(\l^0)=(x_{n},y_{n})$ with initial
point $\lambda^{0}= (x^{0},y^{0})$. Thus $x_{n}\rightarrow 0$ as
$n\rightarrow \infty$. By symmetry of $x$ and $y$ we get
$y_{n}\rightarrow 0$ as $n\rightarrow \infty$.
 \vskip 0.3 truecm

 {\bf Case $q = 1$.} In this case we get $x_{n + 1}\leq x_{n}$, hence
 $$\lim_{n\to \infty}x_n = \alpha \geq 0 \ \ \mbox{exists}.$$
 Similarly,
   $$\lim_{n\to \infty}y_n= \beta\ \ \mbox{also  \ \ exists}.$$
   Thus the point $(\alpha , \beta )$ must be a fixed point for the operator
   (17). Since $\lambda_{0} =(0,0)$ is unique fixed point for $c \leq
 \frac{1}{2}$ (we shall prove uniqueness in section (iii) of this proof), we get $(\alpha,\beta) = (0,0)$.
\vskip 0.3 truecm

 {\bf Remark.} The argument used in the case $q=1$ also works for
 the case $q<1$. But in the case $q < 1$ we proved that the rate of
 convergence to $\lambda_{0}$ is faster than $q^{n}$.

(ii) Invariance of $M_{0},M_{=},M_{1}$ are
 straightforward. Invariance of $M_{<}$, $M_{>}$ follow from the
 following equality
 $$ x'-y' = (x-y)(2cz+ a(x+y)),\ \ \mbox{where}\ \ z=1-x-y\geq 0$$
 which can be obtained from (17).

(iii) Clearly $\lambda_0 = (0,0)$ is a fixed point independently on
parameters $a,b,c$. To get other fixed points consider several
cases: \vskip 0.3 truecm

{\bf Case $x=0, y\neq0$:} From the second equation one gets $y
 =\frac{2c-1}{2c-a}$ which is between $0$ and $1$ iff $ c
 >\frac{1}{2}$. Thus $\lambda_{1} = \left(0,\frac{2c-1}{2c-a}\right)$ is a fixed
 point.
\vskip 0.3 truecm

{\bf Case $x\neq 0, y=0$} is similar to the previous case and gives
$\lambda _{2} =\left(\frac{2c-1}{2c-a},0\right)$.
 \vskip 0.3 truecm

{\bf Case $x\neq 0, y\neq0$:} From (17) one gets a system of linear
 equations, which has unique solution $\lambda _{3}=
 \left(\frac{2c-1}{4c-a-2b},\frac{2c-1}{4c-a-2b}\right)$ (for $c>\frac{1}{2}, a\neq 2b)$.
 Note that if $c\leq\frac{1}{2}$ then
 there is only $\lambda_{0}$.

 To check the type of fixed points
consider Jacobian at $\lambda =(x,y)$
$$\textbf{J}(\lambda)=\textbf{J}(x,y) =\left(
\begin{array}{cccccc}
2c +2(a-2c)x + 2(b-c)y & 2(b-c)x \\[2mm]
2(b-c)y & 2c + 2(a-2c)y + 2(b-c)x
\end{array}
\right).\eqno (19) $$

It is easy to see that the eigenvalues $\mu_{1}(\lambda),
\mu_{2}(\lambda)$ of (19) at fixed points are
$$\mu_{1}(\lambda_{2})=\mu_{2}(\lambda_{1})= \mu_{2}(\lambda_{3})=
2(1-c)<1,$$
$$|\mu_{1}(\lambda_{1})|=
|\mu_{2}(\lambda_{2})|=\left|\frac{2c(1-a) +2b(c-1)}{2c-a}\right|=
\left\{\begin{array}{ll} < 1 \ \ \mbox{if} \ \ a>2b \\
> 1 \ \ \mbox{if} \ \ a<2b
\end{array}\right. $$

$$|\mu_{1}(\lambda_{3})|=\left|\frac{4c(b-1) + 2a(c-1)}{a + 2b -4c}\right|=
\left\{\begin{array}{ll} < 1 \ \ \mbox{if} \ \ a<2b \\
> 1 \ \ \mbox{if} \ \ a>2b
\end{array}\right. $$

This completes the proof of (iii).

(iv) For $a=2b$ the operator (17) has the following form
$$ \left\{\begin{array}{ll}
x'= x (2c + 2(b-c)(x+y))\\
y'= y (2c + 2(b-c)(x+y)).
\end{array}\right.\eqno(20)$$

It is easy to see that $\lambda _{0} =(0,0)$ and any point of
$F=\left\{\lambda =(x,y):x+y=\frac{2c-1}{2(b-c)}\right\}$ is fixed
point if $c> \frac{1}{2}$. Invariance of $I_{\nu}$ follows easily
from the following relation $\frac{y'}{x'}=\frac{y}{x}=\nu$. To
check $\lambda^{(n)}\rightarrow \overline{\lambda}_{\nu}$ for
$\lambda^{0}\in I_{\nu}$, consider restriction of operator (20) on
$I_{\nu}$ which is $x'=\varphi(x)=x\big(2c + 2(b-c)(1 +\nu)x\big)$.
The function $\varphi$ has two fixed points $x=0$ and
$\overline{x}=\frac{1-2c}{2(b-c)(1+\nu)}$. The point $x=0$ is
repeller and $\overline{x}$ is attractive independently on $\nu$
since $\varphi'(\overline{x})=2(1-c)<1$ for $c>\frac{1}{2}$. One can
see that $x^{*}\geq \overline{x}$ where $x^{*}$ is the critical
point i.e $\varphi'(x^{*})= 0$. The graphical analysis shows that
$\overline{x}$ is the global attractive point on $I_{\nu}$.

(v) The existence of $\gamma$ follows from Theorem 1. Other
statements of (v) are straightforward. The theorem is proved.

Note that 2-Volterra operator corresponding to (17) has the
following form
$$
\left\{\begin{array}{ll}
x'=x(ax+2by+2cz)\\
y'=y(2bx+ay+2cz)\\
z'=1-2c(x+y)-(a-2c)(x^{2}+y^{2})-4(b-c)xy
\end{array}\right..\eqno {(21)}
$$\\
Using Theorem 2 one can describe the phase portraits of the
trajectories of (21).

\vskip 0.3 truecm

{\bf Remark.} One of the main goal by introducing the notion of
 $\ell$-Volterra operators was to give an example of QSO which has
 more rich dynamics than Volterra QSO. It is well known \cite{Gan1} that
 for Volterra operators (see (6)) if $a_{ij}\neq 0 \ \ (i\neq j)$ then
 for  any non-fixed initial point $\lambda^{0}$ the set
 $\omega(\lambda^{0})$ of all limit points of the trajectory
 $\{\lambda^{(n)}\}$ is subset of the boundary of simplex. But in our
 case Theorem 2 shows that the limit set need not to be subset of the
 boundary of $S^{2}$.

\end{document}